\documentclass[11pt,reqno]{amsart}
\numberwithin{equation}{section}

\newtheorem{theorem}{Theorem}[section]

\newtheorem{remark}[theorem]{Remark}

\def\al{\aligned}
\def\eal{\endaligned}
\def\be{\begin{equation}}
\def\ee{\end{equation}}
\def\lab{\label}

\def\e{\epsilon}

\def\R{{\bf R}}
\def\M{{\bf M}}

\def\al{\aligned}

\def\pa{\partial}
\def\nb{\nabla}

\numberwithin{equation}{section}

\usepackage[backref=page]{hyperref}

\begin{document}

\title[]{A Sharp Li-Yau gradient bound on Compact Manifolds}
\author{Qi S. Zhang}
\address{Department of Mathematics, University of California, Riverside, CA 92521, USA}
\email{qizhang@math.ucr.edu}
\date{MSC2020: 58J35, 53C44}

\begin{abstract}
Let $(\M^n, g)$ be a $n$-dimensional, complete (compact or noncompact) Riemannian manifold whose Ricci curvature is bounded from below by a constant $-K \le 0$. Let $u$ be a positive solution of the heat equation on $\M^n \times (0, \infty)$. The well known Li-Yau gradient bound states that
$$
t \left(\frac{|\nabla u|^2}{u^2} - \alpha\frac{\pa_t u}{u}\right) \leq \frac{n\alpha^2}{2} +
t \frac{n\alpha^2K}{2(\alpha-1)},\quad \forall \alpha>1, t>0.
$$ The bound with $\alpha =1$ is sharp if $K=0$. If $-K < 0$, the bound tends to infinity if $\alpha=1$. In over 30 years, several sharpening of the bounds have been obtained with $\alpha$ replaced by several functions $\alpha=\alpha(t)>1$ but not equal to $1$. An open question ( \cite{CLN}, \cite{LX} etc) asks if a sharp bound can be reached. In this short note, we observe that for all complete compact manifolds one can take $\alpha=1$. Thus a sharp bound, up to computable constants, is found in the compact case. This result also seems to sharpen Theorem 1.4 in \cite{LY} for compact manifolds with convex boundaries. In the noncompact case one can not take $\alpha=1$ even for the hyperbolic space. An example is also given, which shows that there does not exist an optimal function $\alpha=\alpha(t, K)$ for all noncompact manifolds with negative Ricci lower bound $-K$, giving a negative answer to the open question in the noncompact case.
\end{abstract}
\maketitle

\section{Introduction}

In \cite{LY}, Li and Yau proved the important Li-Yau gradient bound for positive solutions of the heat equation on $\M$ when the Ricci curvature of $\M$ is bounded from below. Here $\M$ is a complete Riemannian manifold of dimension $n$, which may be compact or noncompact. It states that if $Ric\geq -K$ for some constant $K\geq 0$, then for any positive solution $u$ of the heat equation $\pa _t u=\Delta u$, one has
\begin{equation}\label{Li-Yau>1}
\frac{|\nabla u|^2}{u^2} - \alpha\frac{u_t}{u} \leq \frac{n\alpha^2K}{2(\alpha-1)}+\frac{n\alpha^2}{2t},\quad \forall \alpha>1, \, t>0.
\end{equation}
In particular, when $Ric\geq 0$, one obtains the optimal Li-Yau bound
\begin{equation}\label{Li-Yau=1}
\frac{|\nabla u|^2}{u^2} - \frac{u_t}{u} \leq \frac{n}{2t}.
\end{equation}

The Li-Yau gradient bound is a versatile  tool for studying analytical, topological and geometrical properties of manifolds. For instance, the classical parabolic Harnack inequality, optimal Gaussian estimates of the heat kernel, estimates of eigenvalues of the Laplace operator, estimates of the Green's function, and even Laplacian comparison theorem can be deduced from (\ref{Li-Yau=1}).

While the coefficient in \eqref{Li-Yau=1} is sharp for the case where $Ric\geq 0$, the coefficients in \eqref{Li-Yau>1} are not sharp if $K >0$. In the past 30 years, many efforts have been made to improve \eqref{Li-Yau>1}. The readers  may refer to \cite{Ha}, \cite{Ya}, \cite{CTZ}, \cite{Dav}, \cite{GM}, \cite{LX}, \cite{QZZ}, \cite{Wan}, \cite{WanJ}, \cite{BBG}, \cite{YZ}, \cite{Qb} and the references therein for more information. In some of these improvements, the constant $\alpha$ is replaced by  functions in the form of $\alpha(t, K)$ which is strictly greater than $1$ but converges to $1$ as $t \to 0$.
On the other hand, generalizations of Li-Yau gradient bounds have also been studied by many mathematicians. Hamilton \cite{Ha} discovered a matrix Li-Yau type bound for the heat equation. See also \cite{CH} and \cite{CaNi}. Moreover, Li-Yau type bounds were also proved for weighted manifolds with Bakry-\'Emery Ricci curvature being bounded from below, or more generally for metric measure spaces $(X, d, \mu)$ satisfying $RCD^*(K,N)$ condition (see e.g. \cite{BL} and \cite{ZhZx}).

An open question mentioned in several a fore cited papers ( e.g. \cite{CLN} Problem 10.5,  p377;  p4458 \cite{LX} and \cite{BBG} ) asks if a sharp bound can be reached. In this short note, we show that for all complete compact manifolds one can take $\alpha=1$.
 Thus a sharp bound, up to computable constants, is found in the compact case. We point out that going from $\alpha(t, K)>1$ to $\alpha =1$ is a little more than a technical improvement since the latter gives rise to the second order  Laplace estimate for log of solutions.

 Our theorem is

\begin{theorem}\label{main thm}
Let $(\M, g_{ij})$ be a complete, compact $n$-dimensional Riemannian manifold and $u$ a positive solution of the heat equation on $\M \times (0, \infty)$, i.e.,
\be\label{HE}
(\Delta - \pa_t)u=0.
\ee Let $diam_\M$ be the diameter of $\M$ and suppose the Ricci curvature is bounded from below by a non-positive constant  $-K$, i.e. $R_{ij} \ge -K g_{ij}$. Then there exist  dimensional constants $C_1>0$ and $C_2>0$ such that
\be
\lab{sly}
\al
-t \Delta \ln u &= t \left(\frac{|\nabla u|^2}{u^2} - \frac{\pa_t u}{u} \right)\\
&\le \left( \frac{n}{2} + \sqrt{2n K (1+K t)(1+t)} \,  diam_\M + \sqrt{ K (1+K t) (C_1  + C_2 K ) t} \right).
\eal
\ee
\end{theorem}

{\remark The constants $C_1$ and $C_2$ arise from the standard volume comparison or/and heat kernel upper and lower bound which can be efficiently estimated. For compact manifolds, the main utility of the result is for small or finite time. For example, when we choose $C_1$ and $C_2$ to be independent of the volume comparison theorem, the bound actually gives rise to a Laplace and hence volume comparison theorem as $t \to 0$.  For large time, a solution converges to a constant. As $t \to \infty$, it is expected that the bound will be of order $o(t)$ instead of $O(t)$ as of now. But we will not pursue the improvement this time. It would still be interesting to find a sharp form of the RHS of \eqref{sly}.
}

The idea of the proof is the following. First we prove the sharp bound for heat kernels $G$ instead of general positive solutions. This involves an iterated integral estimate on $-\Delta \ln G$ which satisfies a nonlinear heat type equation and Hamilton's estimate of $|\nb \ln G|$ by $t \ln A/G$, taking advantage of the finiteness of the diameter of $\M$ and lower bound of $G$. Here $A$ is a constant which dominates the heat kernel in a suitable time interval. Next we use crucially the result in \cite{YZ} which states that a Li-Yau type bound on the heat kernel implies the same bound on all positive solutions.

 For those who is curious about the noncompact case, we give a negative answer on the existence of a sharp function which depends on time $t$ and Ricci lower bound only and for which the Li-Yau bound holds. First, in the noncompact case one can not prove a bound like \eqref{sly} with a finite right hand side for a fixed $t$ even for the hyperbolic space. For example in the 3 dimensional hyperbolic space with the standard metric, the heat kernel (c.f. \cite{Dav} e.g.) is:
 \[
 G(x, t, y) = \frac{1}{(4 \pi t)^{3/2}} \frac{  d(x, y)}{\sinh d(x, y)} e^{-t- \frac{d^2(x, y)}{4t}}.
 \]It is easy to check, for $r=d(x, y)$,
 \[
 -t \Delta \ln G(x, t, y) = - t \sinh^{-2} r \frac{\partial}{\partial r} \left(\sinh^2 r \frac{\partial}{\partial r} \right) \ln \left(\frac{1}{(4 \pi t)^{3/2}} \frac{r}{\sinh r} e^{-t- \frac{r^2}{4t}}\right)
 \]is of the order $d(x, y)$ when $d(x, y) \to \infty$. See also \cite{YZ2}. Therefore one can not take $\alpha=1$ for all noncompact manifolds with Ricci bounded from below.
 Second, let $\M_0$ be any compact manifold of dimension $n$ and $\M=\M_0 \times \R^1$ be the product manifold of $\M_0$ with $\R^1$. Then the heat kernel on $\M$ is given by $G=G_0(x, t, y) G_1(z, t, w)$ where $G_0$ and $G_1$ are the heat kernel on $\M_0$ and $\R^1$ respectively with $x, y \in \M_0$ and $z, w \in \R^1$ and $t>0$. Fixing $y$ and $w$, we compute
 \[
 \al
 -t \Delta \ln G &= -t (\Delta_0 + \partial^2_z) \ln G = -t\Delta_0 \ln G_0 -  t \partial^2_z \ln G_1\\
 &\le \left( \frac{n+1}{2} + \sqrt{2n K (1+K t)(1+t)} \,  diam_{\M_0} + \sqrt{ K (1+K t) (C_1  + C_2 K ) t} \right).
 \eal
 \]Here we just used Theorem \ref{main thm} on $\M_0$ with $-K$ being the Ricci lower bound and $\Delta_0$ being the Laplace-Beltrami operator on $\M_0$. So,  according to \cite{YZ},  for the product manifold $\M$, the above Li-Yau bound with $\alpha=1$ actually holds for all positive solutions of the heat equation,  which is optimal.  These two examples show that there does not exist a single optimal function of time $t$ and $K$ only, which works for all noncompact manifolds with Ricci curvature bounded from below by a negative constant $-K$. The reason is that for each negative constant $-K$, the optimal function for the above manifold $\M_0 \times \R^1$ is $\alpha=1$. But for the hyperbolic space with $Ricci = - K g$, the optimal function is worse than $\alpha=1$.

We will see at the end of the paper  that a relatively sharp Li-Yau type bound involving $d(x, y)$ is possible for the heat kernels in the noncompact case.

\section{Proof of the theorem}
Let $(\M, g_{ij})$ be a complete, compact Riemannian manifold. In the following, we use $B(x,r)$ and $|B(x,r)|$ to denote the geodesic ball with radius $r$ in $\M$ centered at $x$ and its volume, respectively. The distance between $x, y \in \M$ is denoted by $d(x, y)$.

Now we are ready to prove  Theorem 1.1.

\proof

Let $u=u(x, t)$ be a positive solution of the heat equation on $\M \times (0, \infty)$. Note that no assumption is made on  $u$ at $t=0$. In particular $u$ can be chosen as the heat kernel $G=G(x, t, y)$, $t>0$.
Write
\[
Y=Y(x,t) \equiv \frac{|\nb u|^2}{u^2}-\frac{\pa_t u}{u}= - \Delta \ln u.
\]
According to the computation in Sec.1 of \cite{LY}, we have
\be
\lab{ineqtQ}
\al
(\Delta  - \pa_t) Y
+ 2 \frac{\nb u}{u} \nb Y
&=  2 | \, Hess \, \ln u \, |^2 + 2 Ric ( \nb \ln u, \nb \ln u )\\
& \ge \frac{2}{n} Y^2 - 2 K | \nb \ln u |^2.
\eal
\ee Let $Y^+ \equiv \max \{ Y(x, t), 0 \}$. The inequality \eqref{ineqtQ} implies that $Y^+$ is a subsolution of the inequality in the weak sense: on $\M \times (0, \infty)$,
\be
\lab{ineqtY}
(\Delta  - \pa_t) Y^+
+ 2 \nb \ln u \nb Y^+
 \ge \frac{2}{n} (Y^+)^2 - 2 K | \nb \ln u |^2.
\ee For a positive integer $j$ and a small positive number $\e$, the function $((t-\e)^+)^{2j+2} (Y^+)^{2 j}$ is a legal test function for \eqref{ineqtY} on $\M \times (0, T]$ for any $T>0$. Therefore
\be
\lab{t123}
\al
\frac{2}{n} \int^T_0 \int_{\M}& [(t-\e)^+ Y^+]^{2j+2} dxdt \\
&\le \int^T_0 \int_{\M}  ((t-\e)^+)^{2j+2} (Y^+)^{2 j} (\Delta - \pa_t) Y^+ dxdt \\
&\qquad + \int^T_0 \int_{\M} ((t-\e)^+)^{2j+2} (Y^+)^{2 j} 2 \nb \ln u \nb Y^+ dxdt\\
&\qquad + 2 K \int^T_0 \int_{\M} |\nb \ln u |^2 \, ((t-\e)^+)^{2j+2} (Y^+)^{2 j} dxdt\\
&\equiv T_1 + T_2 + T_3.
\eal
\ee  Let us bound $T_1, T_2$ and $T_3$ respectively.  Using integration by parts, we see that
\be
\lab{t1}
\al
T_1 &= -  2j \int^T_0 \int_{\M}  ((t-\e)^+)^{2j+2} (Y^+)^{2 j-1} |\nb  Y^+|^2 dxdt
- \frac{1}{2j+1} \int_{\M}  ((T-\e)^+)^{2j+2} (Y^+)^{2 j+1}(x, T) dx \\
&\qquad + \frac{2j+2}{2j+1} \int^T_0 \int_{\M} ((t-\e)^+ Y^+)^{2 j+1}  dxdt.
\eal
\ee Writing $(Y^+)^{2j} \nb Y^+ = \frac{1}{2j+1} \nb (Y^+)^{2j+1}$ and doing integration by parts, we deduce
\be
\lab{t2}
\al
T_2 &= \frac{2}{2j+1} \int^T_0 \int_{\M} \nb \ln u \nb (Y^+)^{2j+1}  (t-\e)^+)^{2j+2} dxdt \\
&= - \frac{2}{2j+1} \int^T_0 \int_{\M} (\Delta \ln u ) \, (Y^+)^{2j+1}  (t-\e)^+)^{2j+2} dxdt \\
&= \frac{2}{2j+1} \int^T_0 \int_{\M}  ((t-\e)^+  Y^+)^{2j+2}   dxdt.
\eal
\ee Substituting \eqref{t1} and \eqref{t2} into \eqref{t123} and throwing away the non-positive terms on the right hand side, we find
\be
\lab{a2j+2}
\al
&\left( \frac{2}{n}- \frac{2}{2j+1} \right) \int^T_0 \int_{\M}  ((t-\e)^+  Y^+)^{2j+2}   dxdt\\
&\le \frac{2j+2}{2j+1} \int^T_0 \int_{\M}  ((t-\e)^+  Y^+)^{2j+1}   dxdt
+\underbrace{2 K \int^T_0 \int_{\M} |\nb \ln u |^2 \, ((t-\e)^+)^{2j+2} (Y^+)^{2 j} dxdt}_{T3}.
\eal
\ee

This estimate holds for all positive solutions. Now we take, in particular $u=G(x, t, y)$ the heat kernel
with pole at  a fixed $y \in \M$. We will find a bound for $T_3$ next, which relies on the following estimate by Hamilton \cite{Ha} Theorem 1.1. Let $f$ be a positive solution of the heat equation on $\M \times (0, \infty)$ with $f \le A$ for some constant $A$. Then
\be
\lab{ham1}
t |\nb \ln f |^2 \le (1 + 2 K t) \ln (A/f).
\ee  For a time $t_0>0$ we take $f=f(x, t)=G(x, t+t_0, y)$ with $t \in [0, t_0]$.
According to \cite{LY}, the following upper and lower bound for $f$ hold. See also \cite{Lib} Chapter 13 for further refinement. For some positive constants $C_1, C_2$ depending at worst on the dimension $n$,
\be
\lab{gaulb}
\frac{ \exp (-C_2 K (t+t_0) - \frac{d^2(x, y)}{ 3 (t+t_0)} )}{C_1 |B(x, \sqrt{t+t_0})|} \le f(x, t) \le \frac{C_1 \exp (C_2 K (t+t_0) - \frac{d^2(x, y)}{ 5 (t+t_0)} )}{|B(x, \sqrt{t+t_0})|}.
\ee Note the center of the ball $x$ may be replaced by $y$ in the above bounds by modifying the constants appropriately.
This upper bound implies
\[
A = \sup_{\M \times (0, t_0)} f  \le \frac{C_1 \exp (2 C_2 K t_0 )}{\inf_{z \in \M} |B(z, \sqrt{t_0})|},
\]which yields, by the lower bound of $f$, that
\be
\lab{a/f1}
\frac{A}{f(x, t)} \le \frac{C^2_1 \exp (4 C_2 K t_0 ) \, \exp(\frac{d^2(x, y)}{ 3 t_0} )}{\inf_{z \in \M} |B(z, \sqrt{t_0})|}  |B(x, \sqrt{2 t_0})|.
\ee We mention that $\inf_{z \in \M} |B(z, \sqrt{t_0})|$ can also be replaced by $|B(y, \sqrt{t_0})|$ after adjusting the constants. According to the standard volume comparison theorem, there exists a dimensional constant $C_3$ such that
\[
\frac{|B(x, \sqrt{2 t_0})|}{\inf_{z \in \M} |B(z, \sqrt{t_0})|} \le \exp ( C_3 \sqrt{K} \frac{ diam_{\M} }{\sqrt{t_0}}),
\]where $diam_{\M}$ is the diameter of $\M$. Note that one can also use the heat kernel bounds to prove a similar bound, avoiding the use of volume comparison.
 This and \eqref{a/f1} infer that
\[
\ln \frac{A}{f(x, t)} \le 2 \ln C_1 + 4 C_2 K t_0 + \frac{d^2(x, y)}{ 3 t_0} + C_3 \sqrt{K} \frac{ diam_{\M} }{\sqrt{t_0}},
\]which implies, via \eqref{ham1}, that
\[
t |\nb \ln f(x, t) |^2 \le (1 + 2 K t_0) \left( 4 C_2 K t_0 + 4 C^2_3 K + \frac{diam^2_\M}{ 2 t_0} + 2 \ln C_1 \right), \qquad t \in (0, t_0].
\]Recall that $f(x, t)=G(x, t+t_0, y)$. We can take $t=t_0$ and use the arbitrariness of $t_0$ to conclude
\be
\lab{nblnG}
t |\nb_x \ln G(x, t, y) |^2 \le 2 (1 +  K t) \left( 2 C_2 K t + 4 C^2_3 K + \frac{diam^2_\M}{ t} + 2 \ln C_1 \right), \,  \forall t>0.
\ee

The bound is adequate for us when the time is short, say $t \le 4$.  When $t$ is large, the heat kernel converges to the positive constant $1/|\M|$ where $|\M|$ is the volume of $\M$. In this case the above bound becomes inaccurate. Instead we will use a better bound based on the Li-Yau Harnack inequality (Theorem 2.2 (ii)) in \cite{LY}.  Pick any time $t \ge 4$. Since $\int_\M G(x, t+1, y) dx =1$, there is a point $x_1 \in \M$ such that $G(x_1, t+1, y)=1/|\M|$. According to the afore mentioned Theorem 2.2 (ii) with $q=0$, $\alpha=2$, $t_2=t+1, t_1=t$ therein,
there exists a dimensional constant $C_0>0$ such that
\be
\lab{g1.25}
G(x, t, y) \le G(x_1, t+1, y) ((t+1)/t)^{n}  \exp(C_0 K + \frac{d^2(x, x_1)}{2}).
\ee Since $t \ge 4$, this implies
\[
G(x, t, y) \le 1.25^n |\M|^{-1}  \exp(C_0 K + diam^2_\M/2).
\]Similarly, there is a point $x_2$ such that $G(x_2, t-1, y)=1/|\M|$ and that
\[
G(x_2, t-1, y) \le G(x, t, y) (t/(t-1))^{n}  \exp(C_0 K + \frac{d^2(x, x_2)}{2}),
\]which infers
\be
\lab{g3/4}
G(x, t, y) \ge (3/4)^n |\M|^{-1}  \exp(-C_0 K - diam^2_\M/2), \quad t \ge 4.
\ee Using \eqref{g1.25} and \eqref{g3/4} and recalling $f(x, t)=G(x, t+t_0, y)$, we find, for $t_0 \ge 4$, that
\[
\ln \frac{A}{f(x, t)} \le n \ln 2 + 2 C_0 K + diam^2_\M, \quad t >0.
\] This and \eqref{ham1} yield, for $t \in (0, t_0]$,
\[
t |\nb \ln f(x, t) |^2 \le (1 + 2 K t_0) (n \ln 2 + 2 C_0 K + diam^2_\M).
\]Therefore, after taking $t=t_0$ and renaming, we deduce
\be
\lab{nblnG>4}
t |\nb_x \ln G(x, t, y) |^2 \le 2 (1 +  K t) (n \ln 2 + 2 C_0 K + diam^2_\M), \,  t \ge 8.
\ee

Next we plug \eqref{nblnG} for $t<8$ and \eqref{nblnG>4} for $t \ge 8$ into the term $T_3$ in \eqref{a2j+2} with $u=G(x, t, y)$. After renaming some dimensional constants, we obtain
\[
T_3 \le 4 K (1 +  K T) \left(diam^2_\M + T \, diam^2_\M + C_4 T   +  C_5 K T \right) \int^T_0 \int_{\M}  ((t-\e)^+  Y^+)^{2 j} dxdt
\]Then \eqref{a2j+2} becomes
\be
\lab{a2j+22}
\al
&\left( \frac{2}{n}- \frac{2}{2j+1} \right) \int^T_0 \int_{\M}  ((t-\e)^+  Y^+)^{2j+2}   dxdt
\le \frac{2j+2}{2j+1} \int^T_0 \int_{\M}  ((t-\e)^+  Y^+)^{2j+1}   dxdt\\
&\qquad +4 K (1 +  K T) \left( diam^2_\M + T \, diam^2_\M + C_4 T   +  C_5 K T \right) \int^T_0 \int_{\M}  ((t-\e)^+  Y^+)^{2 j} dxdt.
\eal
\ee Using the notation
\[
A_{j, \e} = \left(\int^T_0 \int_{\M}  ((t-\e)^+  Y^+)^j   dxdt \right)^{1/j},
\]we can write \eqref{a2j+22} as
\be
\lab{a2j+23}
\al
&\left( \frac{2}{n}- \frac{2}{2j+1} \right) A^{2j+2}_{2j+2, \e}\\
&\le \frac{2j+2}{2j+1} A^{2j+1}_{2j+1, \e} +4 K (1 +  K T) \left(diam^2_\M + T \, diam^2_\M + C_4 T   +  C_5 K T \right) A^{2j}_{2j, \e}.
\eal
\ee By H\"older inequality,
\[
A^{2j+1}_{2j+1, \e} \le A^{2j+1}_{2j+2, \e} \left(\int^T_0\int_\M dxdt \right)^{1/(2j+2)}, \quad
A^{2j}_{2j, \e} \le A^{2j}_{2j+2, \e} \left(\int^T_0\int_\M dxdt \right)^{2/(2j+2)},
\]which imply, together with \eqref{a2j+23}, that
\be
\lab{a2j+24}
\al
&\left( \frac{2}{n}- \frac{2}{2j+1} \right) A^2_{2j+2, \e}
\le \frac{2j+2}{2j+1} A_{2j+2, \e} \left(\int^T_0\int_\M dxdt \right)^{1/(2j+2)}\\
&\quad +4 K (1 +  K T) \left(diam^2_\M + T \, diam^2_\M + C_4 T   +  C_5 K T\right) \left(\int^T_0\int_\M dxdt \right)^{2/(2j+2)}.
\eal
\ee Letting $j \to \infty$, we arrive at
\[
\frac{2}{n} A^2_{\infty, \e}
\le  A_{\infty, \e}  +4 K (1 +  K T) \left(diam^2_\M + T \, diam^2_\M + C_4 T   +  C_5 K T \right),
\]where
\[
A_{\infty, \e} = \sup_{\M \times [0, T]} (t-\e)^+ Y^+.
\]This shows
\[
 \sup_{\M \times [0, T]} (t-\e)^+ Y^+ = A_{\infty, \e} \le \frac{n}{2} + \sqrt{2n K (1+K T)
  (diam^2_\M + T \, diam^2_\M + C_4 T   +  C_5 K T )}.
\] Since $\e>0$ is arbitrary and $Y=-\Delta \ln u = - \Delta_x \ln G(x, t, y)$, we conclude that
\be
\lab{lyG}
- t \Delta_x \ln G(x, t, y) \le \frac{n}{2} + \sqrt{2n K (1+K t)(1+t)} \,  diam_\M + \sqrt{ K (1+K t) (C_4 t + C_5 K t)}, \quad \forall t>0.
\ee

From \eqref{lyG}, a short but nice argument from \cite{YZ} implies that the same  bound actually holds if one replaces the heat kernel by any positive solutions of the heat equation. For completeness we present the proof. Let $u$ be any positive solution of the heat equation on $\M \times (0, \infty)$. Without loss of generality we can assume $u_0=u(x, 0)$ is smooth because we can consider $u(x, t+ \e)$ and let $\e \to 0$ otherwise. Then $u(x, t)=\int_\M G(x, t, y) u_0(y) dy$. We compute
\be
\lab{YuZ}
\al
&u^2 (-  \Delta \ln u)(x, t) =  |\nb u|^2 -  u \pa_t u  = |\nb u|^2 -  (1/2) \pa_t u^2\\
&=\int_\M \int_\M \bigg(\nb_x G(x, t, y) \nb_x G(x, t, z) -
\frac{1}{2} \pa_t G(x, t, y)  G(x, t, z) \\
 &\qquad \qquad \qquad - \frac{1}{2}  G(x, t, y)  \pa_t G(x, t, z) \bigg) u_0(y) u_0(z) dydz\\
&\le \int_\M \int_\M \bigg( \frac{1}{2} \frac{|\nb_x G(x, t, y)|^2}{G(x, t, y)} G(x, t, z) + \frac{1}{2} \frac{|\nb_x G(x, t, z)|^2}{G(x, t, z)} G(x, t, y) \\
& \qquad \qquad \qquad
-\frac{1}{2} \pa_t G(x, t, y)  G(x, t, z)
  - \frac{1}{2}  G(x, t, y)  \pa_t G(x, t, z) \bigg) u_0(y) u_0(z) dydz.
\eal
\ee In the above, we have used the Cauchy-Schwarz inequality. Using \eqref{lyG} in the integrand, we conclude that
\[
-t u^2 \Delta \ln u \le \left( \frac{n}{2} + \sqrt{2n K (1+K t)(1+t)} \,  diam_\M + \sqrt{ K (1+K t) (C_4 t + C_5 K t)} \right) u^2.
\]This completes the proof of the theorem. \qed
\medskip

{\remark  By modifying the proof in a standard way to treat the boundary terms,  one can show that the conclusion of Theorem \ref{main thm} (\eqref{sly} there) with perhaps worse parameters depending on the boundary still holds if $\M$  is replaced by a compact domain $D \subset \M$ with smooth convex boundary and the Neumann boundary condition $\pa_n u=0$ is imposed on $u$.  This result sharpens Theorem 1.4 in \cite{LY} for compact manifolds with convex boundary.

Let us indicate the changes in the proof.  Let $n$ be the exterior normal of $\partial D$. It is well known that $\pa_n |\nabla u|^2  = -2 \Pi (\nabla u, \nabla u) \le 0$ on the boundary, where $\Pi$ is the 2nd fundamental form of $\pa D$. Together with the fact that $\pa_n  \ln u=0$ on $\pa D$,  we see that the integral estimate \eqref{a2j+24} is still true since the boundary terms in the steps before \eqref{a2j+24} all can be thrown out.  Hamilton's estimate \eqref{ham1} is also valid for the same reason.  The Neumann heat kernel upper and lower bound in the form of  \eqref{gaulb} in this case is already in \cite{LY} essentially. The upper bound is Theorem 3.2 there. The lower bound is not explicitly stated except when $K=0$ (Theorem 4.2 \cite{LY}). But a rough lower bound like the one in \eqref{gaulb} can be proven similarly or by non-sharp Harnack inequality. Note one needs to use the geodesic balls on $\M$ instead of on $\M \cap D$. The constants in the bound will also depend on $\pa D$.
}

Let us comment on a modification of \eqref{sly} which holds on noncompact manifolds $\M$ with $Ric \ge -K$. Hamilton's estimate \eqref{ham1} is still valid for globally bounded solutions $u>0$ (see \cite{Ko}).  Then assume any version of the Li-Yau bound for the heat kernel $G = G(x, t, y)$ with $\alpha(t, K)>1$, $\alpha(t, K)-1 = K O(t)$ and positive function $\beta=\beta(t, K)$:
\[
t \left(\frac{|\nabla_x G|^2}{G^2} - \alpha(t, K) \frac{\pa_t G}{G}\right) \leq \frac{n \alpha^2(t, K)}{2} +
t \frac{n \beta(t, K) K}{2(\alpha(t, K)-1)},\quad \forall t>0.
\]Using the Gaussian bounds \eqref{gaulb} for $G$ and \eqref{ham1}, one can deduce
\[
\al
t &\left(\frac{|\nabla_x G|^2}{G^2} - \frac{\pa_t G}{G}\right) \leq \frac{n \alpha(t, K)}{2} +
t \frac{n \beta(t, K) \alpha^{-1}(t, K) K}{2(\alpha(t, K)-1)} + C K t^2 \frac{|\nabla G|^2}{G^2}\\
&\le \frac{n \alpha(t, K)}{2} +
t \frac{n \beta(t, K) \alpha^{-1}(t, K) K}{2(\alpha(t, K)-1)} + C_1 K t + C_2 K^2 t^2 + C_3 K d^2(x, y).
\eal
\]However due to the presence of the term $d^2(x, y)$, which is unbounded on noncompact manifolds, this bound can not be translated to all positive solutions. Instead, following \eqref{YuZ}, we obtain, for any positive solution of the heat equation $u=u(x, t)=\int_\M G(x, t, y) u_0(y) dy$, the weaker bound:
\be
\lab{loglapu}
\al
-t \Delta \ln u(x, t) &\le \frac{n \alpha(t, K)}{2} +
t \frac{n \beta(t, K) \alpha^{-1}(t, K) K}{2(\alpha(t, K)-1)} + C_1 K t + C_2 K^2 t^2 \\
&\qquad + \frac{ C_3 K \int_\M G(x, t, y) u_0(y) d^2(x, y) dy}{ \int_\M G(x, t, y) u_0(y) dy}.
\eal
\ee Note this method can also yield a bound similar to \eqref{sly} in the compact case. However the power on the diameter is not sharp. We end the note by remarking that all the dimensional constants can be estimated and improved and the $C_2 K^2 t^2$ term in \eqref{loglapu} can be improved to $C_2 K^2 t$ by mimicking the proof of \eqref{nblnG>4}. If $u_0$ is supported in a ball $B(0, R)$ with $0 \in \M$ and $R>0$, \eqref{loglapu} yields a workable bound, depending on $R$ but not the particular $u_0$.
\be
\lab{loglapu2}
\al
-t \Delta \ln u(x, t) &\le \frac{n \alpha(t, K)}{2} +
t \frac{n \beta(t, K) \alpha^{-1}(t, K) K}{2(\alpha(t, K)-1)} + C_1 K t + C_2 K^2 t^2 \\
&\qquad +  2 C_3 K (d^2(x, 0) + R^2).
\eal
\ee

\medskip

\hspace{-.5cm}{\bf Acknowledgements}
The author gratefully acknowledges the support of Simons'
Foundation grant 710364. Thanks also go to Professors Junfang Li, Peng Lu, Guofang Wei and Meng Zhu for helpful comments.

\end{document}